\documentclass[12pt]{amsart}
\usepackage{geometry}                
\usepackage{graphicx}
\usepackage{amssymb}
\usepackage{epstopdf}
\newtheorem{theorem}{Theorem}
\newtheorem*{slemma}{Stability Lemma}

\title{Systolic inequalities and minimal hypersurfaces}
\author{Larry Guth}
\address{Department of Mathematics, University of Toronto, 40 St.
George St., Toronto ON, Canada}
\email{lguth@math.toronto.edu}

\begin{document}
\begin{abstract} We give a short proof of the systolic inequality for the n-dimensional torus.  The proof uses minimal hypersurfaces.  It is based on the Schoen-Yau proof that an n-dimensional torus admits no metric of positive scalar curvature.
\end{abstract}

\maketitle

In this paper, we give a short new proof of the systolic inequality for the n-dimensional torus.

\begin{theorem} (Gromov) Let $(T^n, g)$ be a Riemannian metric on the n-dimensional torus.  Let $Sys(T^n, g)$ denote the shortest length of a non-contractible curve in $(T^n, g)$.  This length is bounded in terms of the volume by the formula

$$Sys(T^n, g) \le C(n) Vol(T^n, g)^{1/n}.$$

\end{theorem}

This systolic inequality was proven in the 2-dimensional case by Loewner in the 1940's.  The higher-dimensional cases were proven by Gromov in \cite{G} in 1983.  In that paper, Gromov proved a variety of more general results.  Our proof of Theorem 1 does not reproduce all of Gromov's results, but it is pretty short.

Our method gives a stronger result than Gromov's in one respect:  we prove a lower bound for the volume of a metric ball in $(T^n, g)$ of radius comparable to the systole.

\begin{theorem} (Guth) For each dimension $n$, there is a constant $c_n$ so that the following holds.  Suppose that $(T^n, g)$ is any metric on the torus, and suppose that $R < (1/2) Sys(T^n, g)$.  Then there is a point $x$ in $T^n$ so that the ball of radius $R$ around $x$ has volume at least $c_n R^n$.
\end{theorem}

This Theorem 2 immediately implies Theorem 1.    The proof of Theorem 2 goes by induction on the dimension.  The main idea of the paper is that studying the volumes of balls - instead of the total volume - makes the induction work nicely.

Our approach to Theorem 2 is based on an analogy between estimates for volumes of balls and positive scalar curvature.  This analogy appears in the following conjecture of Gromov \cite{G2}:

\newtheorem*{conj}{Conjecture} 

\begin{conj} (Gromov, 1985)  Suppose that $(T^n, g)$ is any metric on the torus, and suppose that $R < (1/2) Sys (T^n, g)$.  Then there is a point $x \in T^n$ so that the ball of radius $R$ around $x$ has volume at least as large as a Euclidean n-ball of radius $R$.
\end{conj}

This conjecture implies Theorem 2.  On the other hand, by taking the limit as $R \rightarrow 0$, it implies that the n-torus admits no metric of positive scalar curvature.

Schoen and Yau invented the minimal hypersurface approach to positive scalar curvature in \cite{SY1} in 1978.  Using this approach, they proved in \cite{SY2} that the n-dimensional torus admits no metric of positive scalar curvature (for $n \le 7$).  Their approach is based on an estimate about the geometry of a stable minimal hypersurface in a manifold of positive scalar curvature.  Our proof of Theorem 2 is an adaptation of their argument.

In dimensions $n \ge 3$, our argument gives a better constant in Theorem 1 than was known previously.

Gromov's paper \cite{G} gives the following estimate:

$$ Sys (T^n, g) \le 6 (n+1) n^n \sqrt{ (n+1)!} \hskip5pt Vol(T^n, g)^{1/n} . $$

This estimate was improved by Wenger in \cite{W}.

$$ Sys(T^n, g) \le 6 \cdot 27^n (n+1)!  \hskip5pt Vol(T^n, g)^{1/n}.$$

Our argument gives the estimate

$$ Sys(T^n, g) \le 8 n \hskip5pt Vol(T^n, g)^{1/n}. $$

A product of $n$ unit circles is a flat n-torus with volume 1 and systole 1.  Remarkably, a randomly chosen flat n-torus with volume 1 has much larger systole, on the order of  $\sim n^{1/2}$.  
Gromov's conjecture above would imply that the systolic inequality holds in the form $Sys(T^n, g) \le C n^{1/2} Vol(T^n, g)^{1/n}$, for an absolute constant $C$.

Besides tori, our Theorem 2 also holds for real projective space $\mathbb{RP}^n$, and for other spaces with nice cohomology rings (see Theorem 3 below).  Gromov's systolic inequality holds for any aspherical or essential manifolds, which makes it much more general.  Moreover, Gromov proved an estimate for the filling radius of a Riemannian manifold in terms of the volume.  The systole of an aspherical manifold is bounded by six times its filling radius, so the filling radius estimate implies the systolic inequality as a corollary.  But the filling radius estimate gives much more geometric information than just a bound for the systole.  The method in this paper doesn't give any estimate for the filling radius.  A bound for the filling radius in terms of the volumes of balls is proven in \cite{Gu}.

Notation: We write $B(p, R)$ to denote the ball around $p$ of radius $R$.  If $Z$ is a d-dimensional surface, we write $|Z|$ to denote its d-dimensional volume.

\section{Proof of the systolic inequality for tori}

We are going to prove the inequality by induction on $n$.  In order to make the induction work, we need to prove something stronger. 

Suppose that $(M, g)$ is a Riemannian manifold, and that $\alpha$ is a cohomology class in $H^1(M, \mathbb{Z}_2)$.  The length of $\alpha$, denoted $L(\alpha)$, is defined to be the smallest length of a 1-cycle $\gamma$ with $\alpha(\gamma) \not= 0$.

\begin{theorem} There exist constants $\epsilon_n > 0$ so that the following holds.  Let $(M^n, g)$ be a closed Riemannian manifold.  Suppose that $\alpha_i \in H^1(M, \mathbb{Z}_2)$ for $1 \le i \le n$ are cohomology classes obeying the following two conditions:

\begin{itemize}

\item{For each $i$, $L(\alpha_i) > 2R$.}

\item{The cup-product $\alpha_1 \cup ... \cup \alpha_n$ is non-zero in $H^n(M^n, \mathbb{Z}_2)$.}

\end{itemize}

Then, for some point $p \in M^n$, the ball around $p$ of radius $R$ has volume at least  $\epsilon_n R^n$.
\end{theorem}

Theorem 3 immediately implies Theorem 2 and Theorem 1.  It also applies to metrics on $\mathbb{RP}^n$.  We will see later that the constant $\epsilon_n$ may be taken to be $(4n)^{-n}$.  Hence any metric on $(T^n, g)$ with systole 2 must have volume at least $(4n)^{-n}$.  Rescaling the result, we see that $Sys(T^n, g) \le 8n Vol(T^n, g)^{1/n}$.

The proof of Theorem 3 relies on a ``stability" estimate for minimizing or near-minimizing hypersurfaces.
It's convenient to use the following definition of near-minimizing hypersurfaces.  Suppose that $(M^n, g)$ is a Riemannian manifold and that $Z$ is an embedded hypersurface in $M$.  We say that $Z$ is
minimizing up to $\delta$ if any embedded surface $Z'$ homologous to $Z$ obeys

$$ Vol (Z) \le Vol(Z') + \delta.$$

\begin{slemma} Let $(M^n, g)$ be a closed Riemannian manifold.  Let $\alpha \in H^1(M, \mathbb{Z}_2)$.  Suppose that $Z^{n-1} \subset (M^n, g)$ is a smooth embedded surface Poincare dual to $\alpha$ and minimizing up to $\delta$ (among embedded surfaces in its homology class).  Suppose that $R < (1/2) L(\alpha)$.  Let $p$ be any point in $M$.  Then the following inequality holds:

$$| Z \cap B(p, R/2) | \le 2 R^{-1} |B(p, R)| + \delta .$$

\end{slemma}

\begin{proof} Consider the spheres around $p$ of radius $t$: $S(p, t)$.  By the coarea formula, we can choose $R/2 < t < R$ so that $|S(p,t)| \le 2 R^{-1} |B(p,R)|$.  Consider $Z$ as a relative cycle in $(B(p,t), S(p,t))$.  We claim that $Z$ has a vanishing relative homology class.  If not, by Poincare duality, it would have non-zero intersection number with an absolute cycle $\gamma$ in $B(p,t)$.  Next we apply Gromov's curve-factoring lemma to control the geometry of $\gamma$.

\newtheorem*{CFL}{Curve-factoring lemma}

\begin{CFL} (Gromov, see \cite{G0} page 290) Let $\gamma$ be any 1-cycle in $B(p, t)$.  Then $\gamma$ is homologous to a sum $\sum_i \gamma_i$ where each $\gamma_i$ has length at most $2 t + \epsilon$, where $\epsilon > 0$ is as small as we like.
\end{CFL}

\begin{proof} It suffices to prove the lemma when $\gamma$ is a topological circle.  We write $\gamma$ as a sum of segments $\gamma = \sum \sigma_i$ where each $\sigma_i$ has length at most $\epsilon$.  Let $N$ be the number of segments.  We have divided a circle into $N$ segments, and so they have a total of $N$ endpoints along the circle, which we label as $p_i$.  We can
do the labeling so that the endpoints of $\sigma_i$ are $p_i$ and $p_{i+1}$, where $i$ is a member of the cyclic group $\mathbb{Z}_N$.  Next we pick a path $c_i$ from $p_i$ to $p$ with length at most $t$.  We define $\gamma_i$ to be the cycle $- c_i + \sigma_i + c_{i+1}$ of length at most $2 t + \epsilon$.  Then $\sum \gamma_i = \gamma$ on the level of chains.  \end{proof}

So if $Z$ gives a non-trivial relative homology class in $B(p,t)$, then we can find a curve $\gamma_i$ of length less than $2R$ with a non-zero intersection number with $Z$.  Since $Z$ is Poincare dual to $\alpha$, we see that $\alpha ( \gamma_i) \not = 0$ and so $L(\alpha) < 2 R$, contradicting our hypothesis.

So the relative cycle $Z \cap B(p, t)$ bounds a relative chain.  Since the cycle is embedded, the chain must be a sum of the components of $B(p, t) - Z$ with multiplicity taken in $\mathbb{Z}_2$.  So we see that $Z$ is homologous to a cycle $Z'$ formed by cutting out $Z \cap B_p(t)$ and gluing in a portion of $S(p,t)$.  The resulting cycle can be smoothed to become a smooth embedded hypersurface.

Since $Z$ is minimizing up to $\delta$, we have $ |Z \cap B(p, t)| \le |S(p,t)| + \delta.$ \end{proof}

Now we give the proof of Theorem 3.

\begin{proof} We will prove Theorem 3 by induction on $n$.  The case $n=1$ is trivial.  The key idea is to look at a minimal hypersurface in $M$.  More precisely, we consider a hypersurface in the homology class Poincare dual to $\alpha_n$.  Every codimension 1 homology class can be realized by a smooth embedded submanifold.  We let $Z^{n-1}$ be a smooth embedded submanifold in the Poincare dual of $\alpha_n$ which is minimizing up to $\delta$.  (Here $\delta > 0$ is a tiny number that we will choose later.)

We will check that $Z$ obeys the hypotheses of the theorem.  By Poincare duality, the cup product $\alpha_1 \cup ... \cup \alpha_{n-1}$ is non-zero in $H^{n-1}(Z, \mathbb{Z}_2)$.  Also, the length of $\alpha_i$ in $Z$ is at least the length of $\alpha_i$ in $(M,g)$, which is at least $2R$.  Therefore, we can find a point $p \in Z$, so that the ball in $Z$ around $p$ of radius $R/2$ has volume at least $\epsilon_{n-1} (R/2)^{n-1}$.  A fortiori, if we intersect $Z$ with the ball in $M$ around $p$ of radius $R/2$, the intersection has area at least $\epsilon_{n-1} (R/2)^{n-1}$.  In other words, we have the following
inequality:

$$\epsilon_{n-1} (R/2)^{n-1} <  |Z \cap B(p, R/2)|.$$

But the stability estimate tells us that

$$ |Z \cap B(p,R/2)| \le 2 R^{-1} |B(p, R)| + \delta .$$

Since $\delta$ is as small as we like, we can conclude that

$$(99/100) 2^{-n} \epsilon_{n-1} R^n < | B(p, R) |. $$  \end{proof}

Remarks about the constants: Unwinding the proof, we get $\epsilon_n \sim exp( - n^2 )$.  This can be improved by tweaking the argument.  Namely, in the stability lemma, we can choose a good radius $t$ in the range $\beta R < t < R$.  We get the conclusion that $| Z \cap B(p, \beta R) | \le (1 - \beta)^{-1} R^{-1} |B(p, R)| + \delta$.  On the other hand, by induction, the volume of some $| Z \cap B(p, \beta R)|$
is at least $\epsilon_{n-1} \beta^{n-1} R^{n-1}$.  Hence

$$\epsilon_n \ge (1 - \beta) \beta^{n-1} \epsilon_{n-1} . $$

The optimal choice of $\beta$ is $\beta = \frac{n-1}{n}$.  In this case, $(1 - \beta) \beta^{n-1} \ge  \frac{1}{4 n}$.  Hence we get $\epsilon_n \ge (4n)^{-n} $.

For comparison, if $(T^n, g)$ is a flat torus obeying the hypotheses of Theorem 2, then it has Euclidean balls, and each ball of radius $R$ has volume $\omega_n R^n$ with $\omega_n \sim n^{-n/2}$.  According to Gromov's conjecture, the sharp constant in Theorem 3 should be $\epsilon_n = \omega_n$.

\end{document}